\documentclass[twoside, 12pt]{article}
\usepackage[russian, english]{babel}
\usepackage{epsfig}
\usepackage{amssymb,amsmath,amsfonts,soul,amsthm,enumerate}
\usepackage{amsfonts}
\usepackage{amsmath}
\usepackage{graphicx}
\usepackage{amsthm}
\usepackage[cp1251]{inputenc}
\usepackage[T2A]{fontenc}
\usepackage{mathrsfs}
\usepackage[english]{babel}
\newtheorem{theorem}{Theorem}[section]
\newtheorem{lemma}{Lemma}[section]

\newtheorem{corollary}{Corollary}[section]
\newtheorem{example}{Example}
\numberwithin{equation}{section}
 \def\@evenhead{\vbox{\hbox to \textwidth{\thepage\hfil\sl\leftmark\strut}\hrule}}
 \def\@oddhead{\vbox{\hbox to \textwidth{\rightmark\hfill\thepage\strut}\hrule}}

\usepackage{hyperref}
\usepackage{alphalph}
\usepackage[subnum]{cases}
\usepackage{array}

\newtheorem{proposition}[theorem]{Proposition}

\title{\bf Non self-adjoint correct restrictions and extensions with real spectrum}
\author{B.N. Biyarov,\, Z.A. Zakarieva,\, G.K. Abdrasheva}

\begin{document}
 \sloppy
\maketitle

\markboth{\hfill{\footnotesize\rm   B.N.~Biyarov, Z.A.~Zakarieva, G.K.~Abdrasheva}\hfill}
{\hfill{\footnotesize\sl  Non self-adjoint correct restrictions and extensions with real spectrum}\hfill}
\vskip 0.3cm

\vskip 0.7 cm

\noindent {\bf Key words:}  maximal (minimal) operator, correct restriction, correct extension, real spectrum, non self-adjoint operator.

\vskip 0.2cm

\noindent {\bf AMS Mathematics Subject Classification:}  35B25,  47A55.
\vskip 0.2cm

%
\noindent {\bf Abstract.} The work is devoted to the study of the similarity of a correct restriction to some self-adjoint operator in the case when the minimal operator is symmetric. The resulting theorem was applied to the Sturm-Liouville operator and the Laplace operator. It is shown that the spectrum of a non self-adjoint singularly perturbed operator is real and the corresponding system of eigenvectors forms a Riesz basis.

\section{\large Introduction}
\label{sec1} 
Let a linear operator $L$ be given in a Hilbert space $H$. The linear equation
\begin{equation}\label{eq1.1}
Lu=f
\end{equation}
is said to be \textit{correctly solvable} on $R(L)$ if $\|u\|\leq C\|Lu\|$ for all $u\subset D(L)$ (where $C>0$ does not depend on $u$) and  \textit{everywhere solvable} if $R(L)=H$. If \eqref{eq1.1} is simultaneously correct and solvable everywhere, then we say that $L$ is a \textit{correct operator}. A correctly solvable operator $L_0$ is said to be \textit{minimal} if $\overline{R(L_0)}\neq H$. A closed operator $\widehat{L}$ is called a \textit{maximal operator} if $R(\widehat{L})=H$ and Ker$\,\widehat{L}\neq \{0\}$. An operator $A$ is called a \textit{restriction} of an operator $B$ and $B$ is said to be an  \textit{extension} of $A$ if $D(A)\subset D(B)$ and $Au=Bu$ for all $u\in D(A)$.

Note that if a correct restriction $L$ of a maximal operator $\widehat{L}$ is known, then the inverses of all correct restrictions of $\widehat{L}$ have in the form \cite{Kokebaev}
\begin{equation}\label{eq1.2}
L_K^{-1}f=L^{-1}f+Kf,
\end{equation}
where $K$ is an arbitrary bounded linear operator from $H$ into  Ker$\,\widehat{L}$.

Let $L_0$ be some minimal operator, and let $M_0$ be another minimal operator related to $L_0$ by the equation $(L_0u, v)=(u, M_0v)$ for all $u\in D(L_0)$ and $v\in D(M_0)$. Then $\widehat{L}=M_0^*$ and $\widehat{M}=L_0^*$ are maximal operators such that $L_0\subset \widehat{L}$ and $M_0\subset \widehat{M}$. 
A correct restriction $L$ of a maximal operator $\widehat{L}$ such that $L$ is simultaneously a  correct extension of the minimal operator $L_0$ is called a \textit{boundary correct extension}.
The existence of at least one boundary correct extension $L$ was proved by Vishik in {\cite{Vishik}}, that is, $L_0\subset L\subset \widehat{L}$.

The inverse operators to all possible correct restrictions $L_K$ of the maximal operator $\widehat{L}$ have the form \eqref{eq1.2}, then $D(L_K)$ is dense in $H$ if and only if $\mbox{Ker}\, (I+K^*L^*)=\{0\}$.
 All possible correct extensions $M_K$ of $M_0$ have inverses of the form
$$
M_K^{-1}f=(L_K^*)^{-1}f=(L^*)^{-1}f+K^*f,
$$
where $K$ is an arbitrary bounded linear operator in $H$ with $R(K)\subset \mbox{Ker}\, \widehat{L}$ such that 
\[\mbox{Ker}\, (I+K^*L^*)=\{0\}.\]

\begin{lemma}[Hamburger {\cite[p.\,269]{Hamburger}}]\label{LemmaA} 
Let $A$ be a bounded linear transformation in $H$ and $N$ a linear manifold. If we write $A(N) = M$  then
$$
A^{*}(M^{\perp})=N^{\perp}\cap R(A^{*}).
$$
\end{lemma}

\begin{proposition}[{\cite[p.\,1863]{Biyarov}}]\label{TheoremB} 
A correct restrictions $L_K$ of the maximal operator $\widehat{L}$ are correct extensions of the minimal operator $L_0$ if and only if  $R(K)\subset \mbox{Ker}\,\widehat{L}$ and $R(M_0)\subset \mbox{Ker}\, K^*$.
\end{proposition}

The main result of this work is the following.

\begin{theorem}\label{Theorem1.3}
Let $L_0$ be symmetric minimal operator in a Hilbert space $H$, $L$ be self-adjoint correct extension of the $L_0$, and $L_K$ be correct restriction of the maximal operator $\widehat{L} (\widehat{L}=L_0^{*})$.  If
$$
R(K^*)\subset D(L), \quad I+KL\geq0,
$$
and $\quad I+KL$ is invertible, where $L$ and $K$ are the operators in representation \eqref{eq1.2}, then $L_K$ similar to a self-adjoint operator.
\end{theorem}

\begin{corollary}\label{Cor1.4}
If $K$ satisfies the assumtions of Theorem \ref{Theorem1.3}, then the spectrum of $L_K$ is real, that is, $\sigma(L_K)\subset \mathbb{R}$. 
\end{corollary}

\begin{corollary}\label{Cor1.5}
If $K$ satisfies the assumtions of Theorem \ref{Theorem1.3} and $L^{-1}$ is the compact operator,  then the system of eigenvectors of $L_K$ forms a Riesz basis in $H$.
\end{corollary}

\begin{corollary}\label{Cor1.6}
The results of Theorem \ref{Theorem1.3} are also valid if conditions ``$I+KL\geq0$ and $\quad I+KL$ is invertible{}'' 
replase to condition  ``$KL\geq0${}'' .
\end{corollary}

\begin{corollary}\label{Cor1.7}
The results of Theorem \ref{Theorem1.3},  Corollary \ref{Cor1.4}-\ref{Cor1.6} are also valid for the $L_K^*$.
\end{corollary}

\section{\large Preliminaries}
\label{sec2} 
In this section, we present some results for correct restrictions and extensions which
are used in Section \ref{sec3_ProofTh1}.

If $A$ is bounded linear transformation from a complex Hilbert space $H$ into itself, then the numerical range of $A$ is by definition the set      
$$
W(A)=\{(Ax, x): \,x\in H, \; \|x\|=1\}.
$$
It is well known and easy to prove that if $\sigma(A)$ denotes the spectrum of $A$, then
$$
\sigma_p(A)\subset W(A), \;\;\, \sigma(A)\subset \overline{W(A)},
$$
for the point spectrum $\sigma_p(A)$ and the spectrum $\sigma(A)$ of $A$,  where the bar indicates closure.
The numerical range of an unbounded operator $A$ in a Hilbert space $H$ is defined as
$$
W(A)=\{(Ax, x): \,x\in D(A), \; \|x\|=1\},
$$
and similarly to the bounded case, $W(A)$ is convex and satisfies $\sigma_p(A)\subset W(A)$. In general, the conclusion $\sigma(A)\subset \overline{W(A)}$ does not surely hold for unbounded operators $A$ (see \cite{Chen}). 

\begin{theorem}[Theorem 2 in {\cite[p.\,181]{Radjavi}}]\label{The2.3}
The following are equivalent conditions on an operator $T$:

$(1)$ $T$ is similar to a self-adjoint operator.

$(2)$ $T=PA$, where $P$ is positive and invertible and $A$ is self-adjoint.

$(3)$ $S^{-1}TS=T^*$ and $0\neq \overline{W(S)}$.
\end{theorem}

\begin{theorem}[Theorem 1 in {\cite[p.\,215]{Williams}}]\label{The2.4}
Let $A$ and $B$ operators on the complex Hilbert space $H$. If $0\notin \overline{W(A)}$ then
$$
\sigma(A^{-1}B)\subset \overline{W(B)}/\overline{W(A)}.
$$
\end{theorem}

\begin{corollary}[Corollary in {\cite[p.\,218]{Williams}}]\label{Corol2.5}
If $A>0, \;\, B\geq 0$ and $C=C^*$, then $\sigma(AB)$ is positive and $\sigma(AC)$ is real.
\end{corollary}

\begin{theorem}[Theorem A in {\cite[p.\,508]{Mehdi}}]\label{The2.6}
The numerical range $W(T)$ of $T$ is convex and $W(aT + b) = aW(T) + b$ for all complex numbers $a$ and $b$. 
\end{theorem}


\section{\large Proof of Theorem \ref{Theorem1.3}}
\label{sec3_ProofTh1}

We transform \eqref{eq1.2} to the form
\begin{equation}\label{eq3.1}
L_K^{-1}=L^{-1}+K=(I+KL)L^{-1}.
\end{equation}
Then $L_K$ is defined as the restriction of the maximal operator $\widehat{L}$ on the domain
$$
D(L_K)=\{u\in D(\widehat{L}):\, (I-K\widehat{L})u\in D(L)\}.
$$
Now let us prove Theorem \ref{Theorem1.3}. 
It was proved in \cite[p.\,27]{Biyarov2}  that $KL$ is bounded on $D(L)$ (that is, $\overline{KL}\in B(H)$) if and only if 
$$
R(K^*)\subset D(L^*).
$$
It follows from $\overline{D(L)}=H$  that $\overline{KL}$ is bounded on $H$. In the future, instead of $\overline{KL}$, we will write $KL$. Then, by virtue of Theorem \ref{The2.3} and taking into account the conditions of Theorem \ref{Theorem1.3} that $I+KL\geq 0$ and $I+KL$ is invertible, we obtain proof of Theorem \ref{Theorem1.3}.

The proof of Corollary \ref{Cor1.4} follows from Corollary \ref{Corol2.5}.
Corollary \ref{Cor1.5} is easy to obtain from the fact that the operator 
$$
C=(I+KL)^{1/2}L^{-1}(I+KL)^{1/2}
$$
is self-adjoint and
\begin{equation}\label{eq3.2}
L_K^{-1}=(I+KL)^{1/2}C(I+KL)^{-1/2}=(I+KL)L^{-1}.
\end{equation}

Let us proof Corollary \ref{Cor1.6}. By Theorem \ref{The2.6}, we get that $0\notin W(I+KL)$. Then $I+KL\geq 0$ and $I+KL$ is invertible.

The proof of Corollary \ref{Cor1.7} follows from \eqref{eq3.2}, since $C$ is a self-adjoint operator and in the case Corollary \ref{Cor1.5} the self-adjoint operator $C$ is compact. 


\section{\large Non self-adjoint perturbations for some differential operators}
\label{sec4}
\begin{example}\label{Ex1}
{\rm We consider the Sturm-Liouville equation on the interval $(0, 1)$
\begin{equation}\label{eq4.1}
\widehat{L}y=-y''+q(x)y=f,
\end{equation}
where $q(x)$ is the real-valued function of  $L^2(0, 1)$. We denote by $L_0$  the minimal operator and by $\widehat{L}$  the maximal operator generated by the differential equation \eqref{eq4.1} in the space $ L_2 (0,1) $. It's clear that
\[D(L_0)=\textit{\mbox{\r{W}}}_2^2(0,1)\]
and
\[D(\widehat{L})=\{y \in L^2(0, 1): \ y,\ y' \in AC [0, 1], \  y'' - q(x)y \in L^2(0, 1)\}.\]
Then $\mbox{Ker}\, \widehat{L}=\{a_{11}c(x)+a_{12}s(x)\}$, where $a_{11}, \: a_{12}$ are arbitrary constants, and the functions  $c(x)$ and $s(x)$ are defined as follows
\[c(x)=1+\int_0^x \mathscr{K}(x,t; 0)\,dt, \quad s(x)=x+\int_0^x \mathscr{K}(x,t; \infty )t\,dt,\]
where 
$$
\mathscr{K}(x,t; 0)=\mathscr{K}(x,t)+\mathscr{K}(x,-t), \, \;\; \mathscr{K}(x,t; \infty )=\mathscr{K}(x,t)-\mathscr{K}(x,-t),
$$
and $\mathscr{K}(x,t)$ is the solution of the following Goursat problem
$$
\left \{
\begin{split}
& \dfrac{\partial^2 \mathscr{K}(x,t)}{\partial x^2}-\dfrac{\partial^2 \mathscr{K}(x,t)}{\partial t^2}=q(x)\mathscr{K}(x,t), \\
& \mathscr{K}(x,-x)=0, \quad \mathscr{K}(x,x)=\dfrac{1}{2} \int_{0}^{x} q(t)dt, \\
\end{split}
\right.
$$
in the domain
\[ \Omega=\bigl \{(x,t): \; 0<x<1, \; -x<t<x \bigr \}. \]
Note that $c(0)=s'(0)=1, \; c'(0)=s(0)=0$  and Wronskian
\[W(c,s)\equiv c(x)s'(x)-c'(x)s(x)=1.\]

As a fixed boundary correct extension $L$ we take the operator corresponding to the Dirichlet problem for equation \eqref{eq4.1} on $(0,1)$. Then
\[D(L)=\bigl \{y\in W_2^2(0,1): \; y(0)=0, \: y(1)=0 \bigr \}. \]
Therefore the description of the inverse of all correct restrictions $L_K$ of the maximal operator $\widehat{L}$ has the form
\begin{eqnarray*} 
y\equiv L_K^{-1}f&=&\int_0^x \big[c(x)s(t)-s(x)c(t)\big]f(t)\,dt \\[5pt]
& & \quad - \,\frac{s(x)}{s(1)}\int_0^1\big[c(1)s(t)-s(1)c(t)\big]f (t) \,dt\\[5pt]
& &\quad +\, c(x) \int_0^1 f(t)\overline{\sigma_1(t)}dt+s(x)\int_0^1 f(t)\overline{\sigma_2(t)}\,dt,
\end{eqnarray*}
where $\sigma_1(x), \: \sigma_2(x)\in L_2(0,1)$ which uniquely determine the operator $K$ from \eqref{eq1.2} in the following form
$$
Kf=c(x)\int_0^1 f(t)\overline{\sigma_1(t)}dt+s(x)\int_0^1 f(t)\overline{\sigma_2(t)}dt, \quad \mbox{for all}\; f\in L_2(0,1).
$$
$K$ is a bounded operator in $L_2(0,1)$ acting from $L_2(0,1)$ to $\mbox{Ker}\,\widehat{L}$. The operator$L_K$ is the restriction of  $\widehat{L}$ on the domain
\[
\begin{split}
& D(L_K)=\biggl \{y\in W_2^2(0,1): \; y(0)=\int_0^1 \big[-y''(t)+q(t)y(t)\big]\overline{\sigma_1(t)}dt;\\
& \qquad \qquad\qquad\qquad y(1)=c(1)y(0)+s(1) \int_0^1 \big[-y''(t)+q(t)y(t)\big]\overline{\sigma_2(t)}dt \biggr \}.\\
\end{split}
\]}
\end{example}

From the condition
$$
R(K^*)\subset D(L^*)=D(L)
$$
we have that
$$
KLy=c(x)\int_0^1 y(t)[-\overline{\sigma}_1''(t)+q(t)\overline{\sigma}_1(t)] dt+ s(x)\int_0^1 y(t)[-\overline{\sigma}_2''(t)+q(t)\overline{\sigma}_2(t)] dt,
$$
where
$$
y\in D(L), \;\; \sigma_1, \, \sigma_2\in W_2^2(0, 1), \;\; \sigma_1(0)=\sigma_1(1)=\sigma_2(0)=\sigma_2(1)=0.
$$
If $I+KL\geq0$ and $I+KL$ is invertible, then the spectrum of the operator $L_K$ consists only of real  eigenvalues $\{\lambda_k\}_{k=1}^\infty$ and the corresponding eigenfunctions $\{\varphi_k\}_{k=1}^\infty$ forms a Riesz basis in $L^2(0, 1)$, since $L^{-1}$ is a compact self-adjoint positive operator.
In particular, if
$$
\sigma_1(x)=\alpha(L^{-1}c)(x), \;\;\, \sigma_2(x)=\beta(L^{-1}s)(x), \;\; \, \alpha, \; \beta\geq0,
$$
then $KL\geq0$.
Therefore, by Corollary \ref{Cor1.6}, the results of Theorem \ref{Theorem1.3} are valid for $L_K$.
In this case, $L^{-1}_K$ has the form
$$
y=L_K^{-1}f=L^{-1}f+c(x)\int_0^1 f(t)(L^{-1}c)(t) dt+s(x)\int_0^1 f(t)(L^{-1}s)(t) dt.
$$
Then $(L_K^{-1})^*=(L_K^*)^{-1}$ has form
$$
v(x)=(L^{-1}f)(x)+\alpha(L^{-1}c)(x)\int_0^1 f(t)c(t) dt+\beta(L^{-1}s)(x)\int_0^1 f(t)s(t) dt.
$$
Thus, we have
$$
(L_K^*v)(x)=-v''(x)+q(x)v(x)+a(x)v'(0)+b(x)v'(1)=f(x),
$$
$$
D(L_K^*)=\{v\in W_2^2(0, 1):\, v(0)=v(1)=0\},
$$
where
$$
a(x)=\dfrac{\alpha\beta(c, s)s(x)-\alpha(1+\beta\|s\|^2)c(x)}{(1+\alpha\|c\|^2)(1+\beta\|s\|^2)-\alpha\beta|(c, s)|^2},
$$

$$
b(x)=\dfrac{\alpha[c(1)(1+\beta\|s\|^2)-\beta s(1)(s, c)]c(x)-\beta[\alpha c(1)(c, s)-s(1)(1+\alpha\|c\|^2)]s(x)}{(1+\alpha\|c\|^2)(1+\beta\|s\|^2)-\alpha\beta|(c, s)|^2},
$$
$a(x), \, b(x)\in$ Ker $\widehat{L}$ and $(\cdot, \cdot)$  is scalar product in $L^2(0, 1)$.
The operator $L_K^*$ acts as
$$
L_K^*=L^*+Q,
$$
where
$$
L^*=-\frac{d^2}{dx^2}+q(x),
$$
$$
(Qv)(x)=a(x)<\delta'(x), v(x)>+b(x)<\delta'(x-1), v(x)>=a(x)v'(0)+b(x)v'(1),
$$
that is, the function $Q\in W_2^{-2}(0, 1)$. Thus, we have constructed an example of a non self-adjoint singularly perturbed Sturm-Liouville operator with a real spectrum  and  the system of eigenvectors that forms a Riesz basis in $L^2(0, 1)$.
\begin{example}\label{Ex2}
{\rm In the Hilbert space $L^2(\Omega)$, where  is a bounded domain in $\mathbb{R}^m$ with an infinitely smooth boundary $\partial\Omega$, let us consider the minimal $L_0$  and maximal $\widehat{L}$ operators generated by the Laplace operator
\begin{equation}\label{eq4.2}
-\Delta u=-\biggl(\frac{\partial^2 u}{\partial{x_1^2}}+\frac{\partial^2 u}{\partial{x_2^2}}+\cdots+\frac{\partial^2 u}{\partial{x_m^2}}\biggr).
\end{equation}
The closure $L_0$, in the space $L^2(\Omega)$ of  Laplace operator \eqref{eq4.2} with the domain $C_0^\infty (\Omega)$, is \textit{the minimal operator corresponding to the Laplace operator}.
The operator $\widehat{L}$, adjoint to the minimal operator $L_0$ corresponding to  Laplace operator, is \textit{the maximal operator corresponding to the Laplace operator}. Then
 \[D(\widehat{L})=\{u\in L^2(\Omega): \; \widehat{L}u =-\Delta u \in L^2(\Omega)\}.\]
Denote by $L$ the operator, corresponding to the Dirichlet problem with the domain
\[ D(L)=\{u\in W_2^2(\Omega): \; u|_{\partial\Omega}=0\}. \]
We have \eqref{eq1.2}, where $K$ is an arbitrary linear operator bounded in $L^2(\Omega)$ with 
\[R(K)\subset \mbox{Ker}\,\widehat{L} =\{u\in L^2(\Omega): \: -\Delta u=0\}.\]
Then the  operator $L_K$ is defined by
\[\widehat{L}u= -\Delta u,\]
on
\[D(L_K)=\{ u\in D(\widehat{L}) : \; [(I-K\widehat{L})u]|_{\partial\Omega}=0 \},\]
where $I$ is the identity operator in $L^2(\Omega)$. 
Note that $L^{-1}$ is a self-adjoint compact operator. 
If $K$ satisfies the conditions of Theorem \ref{Theorem1.3}, then $L_K$ is non self-adjoint operator with  a real positive spectrum (i.e., $\sigma(L_K)\subset \mathbb{R}_{+}$), and the system of eigenvectors $L_K$ forms a Riesz basis in $L^2(\Omega)$. In particular, if
$$
Kf=\varphi(x)\int_{\Omega}f(t)\psi(t) dt,
$$
where $\varphi\in W_{2, loc}^2(\Omega)\cap L^2(\Omega)$ is a harmonic function and $\psi\in L^2(\Omega)$, then $K\in B(L^2(\Omega))$ and $R(K)\subset \mbox{Ker}\,\widehat{L}$. From $R(K^*)\subset D(L)$ it follows that $\psi\in W_2^2(\Omega)$ and $\psi|_{\partial \Omega}=0$.
From the condition $KL\geq0$ we have that $\psi(x)=\alpha(L^{-1}\varphi)(x), \;\; \alpha\in \mathbb{R_{+}}$.
Hence the operator $L_K$ is the restriction of $\widehat{L}$ to the domain
\[
D(L_K)= \Big\{u\in D(\widehat{L}):\, \Big(u-\frac{\varphi}{1+\|\varphi\|^2}\int_{\Omega} u(y)\varphi(y) dy\Big)\Big|_{\partial\Omega}=0\Big\}.
\]
The inverse of $L_K^{-1}$ has the form
\begin{equation}\label{eq4.3}
u=L_K^{-1}f=L^{-1}f+\varphi\int_{\Omega}f(y)(L^{-1}\varphi)(y) dy.
\end{equation}
We find the adjoint operator $L_K^*$. From \eqref{eq4.3} we have
$$
v=(L_K^{-1})^*g=L^{-1}g+L^{-1}\varphi \int_{\Omega}g(y)\varphi(y) dy, \quad \mbox{for all} \;\; g\in L^2(\Omega).
$$
Then
\begin{eqnarray*} 
& & L_K^*v=-\Delta v+\frac{\varphi}{1+\|\varphi\|^2}\int_{\Omega} (\Delta v)(y)\varphi(y) dy=g, \\[7pt]
& & D(L_K^*)=D(L)=\big\{v\in W_2^2(\Omega):\, v\big|_{\partial \Omega}=0\big\}.
\end{eqnarray*}
By virtue of Corollary \ref{Cor1.7}, the spectrum of the operator $L_K^*$ consists only of real positive eigenvalues and the corresponding eigenfunctions forms a Riesz basis in $L^2(\Omega)$.
Note that
$$
(L_K^*v)(x)=-(\Delta v)(x)+\frac{\varphi(x)}{1+\|\varphi\|^2}F(u)=g(x),
$$
where $F\in W_2^{-2}(\Omega)$, since
$$
F(u)=\int_{\Omega}(\Delta v)(y)\varphi(y) dy.
$$
This is understood in the sense of the definition of the space $H^{-s}(\Omega), \;\, s>0$ as in Theorem 12.1 (see \cite[p.\,71]{Lions})}. 
\end{example}

Thus, we have shown the examples of a non self-adjoint singularly perturbed operator with a real spectrum. Moreover, the corresponding eigenvectors forms a Riesz basis in $L^2(\Omega)$.


\vskip 1 cm \footnotesize
\begin{flushleft}
   Bazarkan Nuroldinovich Biyarov, \\
   Zaruet Almazovna Zakarieva, \\
   Gulnara Kaparovna Abdrasheva, \\
   Faculty of Mechanics and Mathematics\\
   L.N. Gumilyov Eurasian National University \\
   13 Munaitpasov St,\\
   010008 Nur-Sultan, Kazakhstan\\
   E-mail: bbiyarov@gmail.com,\\
   \qquad\quad\,  zaruet.zakarieva@mail.ru,\\
   \qquad\quad\,  gulnara.abdrash@gmail.com
\end{flushleft} 


\begin{thebibliography}{99.}
\label{bib}
\bibitem{Kokebaev} 
B.K.~Kokebaev,  M.~Otelbaev, A.N.~Shynibekov,  \textit{About expansions and restrictions of operators in Banach space.} Uspekhi Matem. Nauk 37  (1982), no. 4, 116--123 (in Russian).

\bibitem{Vishik} 
M.I.~Vishik, \textit{On general boundary problems for elliptic differential equations}. Tr. Mosk. Matem. Obs. 1 (1952), 187--246  (in Russian). English transl.: Am. Math. Soc., Transl., II, 24 (1963), 107--172.

\bibitem{Hamburger} 
H.L.\,Hamburger, \textit{Five notes on a generalization of quasi-nilpotent transformations in Hilbert space}, Proc. London Math. Soc. \textbf{3},  1  (1951), 494--512.

\bibitem{Biyarov}  
B.N.~Biyarov, \textit{Spectral properties of correct restrictions and extensions of the Sturm-Liouville operator}. Differenisial'nye Uravneniya  30 (1994), no. 12, 2027--2032 (in Russian). English transl.: Differ. Equations  30 (1994), no. 12, 1863--1868.

\bibitem{Chen} D.Y.\,Wu and A.\,Chen, 
\textit{Spectral inclusion properties of the numerical range in a space with an indefinite metric}, Linear Algebra Appl. \textbf{435} (2011), 1131--1136.

\bibitem{Radjavi}H.\,Radjavi and J.P.\,Williams, 
\textit{Products of self-adjoint operators}, Michigan Math. J. \textbf{16} (1969), 177--185. 

\bibitem{Williams}J.P.\,Williams, 
\textit{Spectra of products and numerical ranges}, J. Math. Anal. Appl. \textbf{17} (1967), 214--220.

\bibitem{Mehdi}M.\,Radjabalipour and H.\,Radjavi, 
\textit{On the geometry of numerical ranges}, Pacific J. Math. \textbf{61} (1975), No. 2, 507--511.

\bibitem{Biyarov2} B.N. Biyarov, D.A. Svistunov, G.K. Abdrasheva, 
\textit{Correct singular perturbations of the Laplace operator}, Eurasian Math. J. \textbf{11} (2020), No. 4, 25--34.

\bibitem{Lions} J.L. Lions, E. Magenes 
\textit{Non-Homogeneous Boundary Value Problems and Applications I}. Springer-Verlag, Berlin, 1972.




%
\bigskip
%
\end{thebibliography}
\end{document}